\documentclass[12pt,a4paper]{amsart}
\usepackage[cp850]{inputenc}
\usepackage{amssymb}
\usepackage[mathscr]{eucal}
\usepackage{amscd}
\usepackage{amsmath}
\usepackage[all]{xy}
\usepackage{mathptmx}
\usepackage{color}

\newcommand{\C}{\mathbb C}

\newcommand{\R}{\mathbb R}
\newcommand{\N}{\mathbb N}

\newcommand{\KP}{{\sf K}\hspace{-1pt}{\sf P}}

\newcommand{\Dom}{\mathrm{Dom}}
\newcommand{\lop}{\curvearrowright }

\newcommand{\adef}{\begin{defn}}
\newcommand{\zdef}{\end{defn}}

\newcommand{\Ran}{\mathop{\mathrm{Ran}}}

\newcommand{\To}{\longrightarrow}

\theoremstyle{plain}
\newtheorem{theorem}{Theorem}[section]
\newtheorem{prop}[theorem]{Proposition}
\newtheorem{lemma}[theorem]{Lemma}

\newtheorem{defn}[theorem]{Definition}

\theoremstyle{remark}

\newtheorem{quest}{Question}

\def\PO{\operatorname{PO}}
\def\PB{\operatorname{PB}}

\begin{document}

\title{Operators on the Kalton-Peck space $Z_2$}

\author[J.M.F. Castillo]{Jes\'us M. F. Castillo}
\address{Universidad de Extremadura, Instituto de Matem\'aticas Imuex,
E-06011 Badajoz, Spain.}
\email{castillo@unex.es}

\author[M. Gonz\'alez]{Manuel Gonz\'alez}
\address{Departamento de Matem\'aticas,
Universidad de Cantabria, E-39071 Santander, Spain.}
\email{manuel.gonzalez@unican.es}

\author[R. Pino]{Ra\'ul Pino}
\address{Departamento de Matem\'aticas,
Universidad de Extremadura, E-06011 Badajoz, Spain.}
\email{rpino@unex.es}

\thanks{2010 Mathematics Subject Classification: 46B03, 46B10, 46M18, 47A53, 47L20.\\
Keywords: The Kalton-Peck space; Operators on sequence spaces; Operator ideals; Strictly singular operators.\\
The research of the first and third author was supported in part by project IB20038 funded by Junta de Extremadura.
The research of the first two authors was supported in part by
project PID2019-103961GB funded by MINCIN. The research of the third author was partially supported by project
FEDER-UCA18-108415 funded by 2014-2020 ERDF Operational Programme and by the Department of Economy, Knowledge,
Business and University of the Regional Government of Andalucia.}
\date{\today}

\begin{abstract} We study operators on the Kalton-Peck Banach space $Z_2$ from various points of view: matrix representations, examples, spectral properties and operator ideals. For example, we prove that there are non-compact, strictly singular operators acting on $Z_2$, but the product of two of them is a compact operator. Among applications, we show that every copy of $Z_2$ in $Z_2$ is complemented, and each semi-Fredholm operator on $Z_2$ has complemented kernel and range, the space $Z_2$ is $Z_2$-automorphic and we give a partial solution to a problem of Johnson, Lindenstrauss and Schetchman about strictly singular perturbations of operators on $Z_2$.
\end{abstract}

\maketitle
\thispagestyle{empty}

\section{Introduction}
In spite of being a by now classical Banach space, many aspects of the celebrated Kalton-Peck space $Z_2$ introduced in \cite{kaltpeck} remain unknown. For instance, it is not known the answer to the \emph{hyperplane problem}: is $Z_2$ is isomorphic to its closed subspaces of codimension $1$? The space is clearly  isomorphic to its closed subspaces of codimension $2$. For a  short, sharp, shocked exposition of the few known facts about the structure of $Z_2$ and its subspaces we refer to \cite[Section 10.8]{hmbst} (The Kalton-Peck spaces) and \cite[Section 10.9]{hmbst} (Properties of $Z_2$ explained by itself) and, for a condensed exposition, to the remainder of this section.\medskip

In this paper we study operators on the  space $Z_2$ from several points of view. In Section \ref{sect:matrix} we introduce a  matrix representation of these operators, and the results obtained are refined for operators that admit an upper or lower triangular representation in Section \ref{sect:triangular}. In Section \ref{sect:involution} we consider an involution $T\To T^+$ on the algebra $\mathfrak L(Z_2)$ of bounded operators on $Z_2$ introduced in \cite{kaltsym} and obtain several new results that we apply to show that every copy of $Z_2$ in $Z_2$ is complemented, that the space $Z_2$ is $Z_2$-automorphic and that each semi-Fredholm operator on $Z_2$ has complemented kernel and range. Section \ref{sect:examples} contains a list of natural examples of operators on $Z_2$. In Section \ref{sect:ideals} we follow ideas of Kalton to show that the operator ideals of strictly singular, strictly cosingular and inessential operators coincide and contain each proper  operator ideal of $\mathfrak L(Z_2)$,  we prove that the product of two of strictly singular operators on $Z_2$ is  compact, although there are non-compact, strictly singular operators, and we show that the perturbation classes problem has a positive solution for operators on $Z_2$. The final Section \ref{sect:problems} outlines a few open directions of research.\medskip

Next we summarize some facts about the space $Z_2$ and operators on it. We refer the reader to the Preliminaries (Section \ref{sect:prelim}) for any unexplained notation.

\subsection{Facts about the space $Z_2$} Let us make a brief description of the space $Z_2$.
Consider the homogeneous map $\KP: \ell_2 \to \ell_\infty$ given by $\KP(x) =  2 x \log |x|$ for normalized $x\in \ell_2$. Then
\begin{equation*} Z_2 = \{(\omega, x)\in \ell_\infty\times \ell_2: \omega - \KP x \in \ell_2\}\end{equation*} endowed with the quasinorm $\|(\omega, x)\|_{Z_2} = \|\omega - \KP x\|_2 + \|x\|_2$. The space $Z_2$ is a nontrivial \emph{twisted Hilbert} space in the sense that it contains an uncomplemented subspace $M$ isomorphic to $\ell_2$ such that $Z_2/M$ is again isomorphic to $\ell_2$. Indeed, there is a nontrivial exact sequence
\begin{equation}\label{ip-seq}
\begin{CD} 0@>>>  \ell_2@>i>> Z_2@>p>> \ell_2@>>>0,\end{CD}
\end{equation}
with inclusion $i y= (y,0)$ and quotient map $p(\omega, x)=x$. The quasinorm above is equivalent to a norm \cite[Theorem 4.7]{kaltpeck} and thus $Z_2$ is a reflexive and $\ell_2$-saturated Banach space (see \cite{CastilloGon97}) isomorphic (actually isometric, see below) to its dual \cite[Theorem 5.1]{kaltpeck}, something we will write as $Z_2 \simeq Z_2^*$. We will also consider the following quasi-Banach spaces:\medskip
\begin{itemize}
\item $\Dom \KP= \{x\in \ell_2: \KP x \in \ell_2\}$ endowed with the quasinorm $\|x\| = \|\KP x\|_2 + \|x\|_2$. This quasinorm is equivalent to the norm of the Orlicz space $\ell_f$ generated by the Orlicz function $f(t)= t^2 \log^2 t$ \cite{kaltpeck}; so we can identify $\Dom \KP=\ell_f$.
\item $\Ran \KP= \{\omega\in \ell_\infty: \exists x\in \ell_2, \omega - \KP x\in \ell_2\}$ endowed with the quasinorm $\|\omega\|_{f^*}= \inf_{x\in \ell_2} \|\omega - \KP x\|_2 + \|x\|_2$, and we can identify $\Ran \KP= \ell_f^*$ \cite{cabenon}.\medskip
\end{itemize}

There is another natural nontrivial exact sequence
\begin{equation}\label{jq-seq}
\begin{CD} 0@>>> \ell_f@>j>> Z_2@>q>> \ell_f^* @>>>0,\end{CD}
\end{equation}
with inclusion $jx= (0,x)$, quotient map $q(\omega, x)=\omega$ and associated quasilinear map $\KP^{-1}: \ell_f^* \To\ell_2$ which provides a description of $Z_2$ as a twisted sum of $\ell_f$ and $\ell_f^*$:

 \begin{equation*} Z_2 = \{(x, \omega)\in \ell_\infty\times \ell_f^*: x - \KP^{-1} \omega \in \ell_f\}\end{equation*} endowed with the quasinorm  $\|(x, \omega)\|= \|x-\KP^{-1} \omega\|_{\ell_f} +\|\omega\|_{\ell_f^*}$, which is equivalent to the original one on $Z_2$ and thus there exists $m,M>0$ so that
\begin{equation}\label{equiv-qnorms}
m \| (\omega, x) )\|_{Z_2}  \leq  \|x-\KP^{-1} \omega\|_{\ell_f} +\|\omega\|_{\ell_f^*} \leq
M \| (\omega, x) )\|_{Z_2} \end{equation}
In this representation, $\omega \to (\omega, \KP^{-1} \omega)$ is a bounded homogeneous lifting for $q$ (see \cite{symmetries21}) and\medskip

\begin{itemize}
\item $\Dom \KP^{-1}= \{\omega\in \ell_f^*: \KP^{-1} \omega \in \ell_f\}$ is endowed with the quasinorm $\|\omega\| = \|\KP^{-1} \omega\|_{\ell_f} + \|\omega\|_{\ell_f^*}$. This quasinorm is equivalent to the norm of $\ell_2$ so we can identify $\Dom \KP^{-1}=\ell_2$.
\item $\Ran \KP^{-1}= \{x \in \ell_\infty: \exists \omega \in \ell_f^*, x - \KP^{-1} \omega \in \ell_f\}$ endowed with the quasinorm $\|x\| = \inf_{\omega\in \ell_f^*} \|x - \KP^{-1} \omega\|_{\ell_f} + \|\omega\|_{\ell_f^*}$, and one can identify $\Ran \KP = \ell_2$.\medskip
\end{itemize}

There is no known explicit formula for $\KP^{-1}$ (see \cite{racsam}). Both compositions $\KP \circ \KP^{-1}: \ell_f^* \to \ell_f^*$ and $\KP^{-1}\circ \KP: \ell_2\to \ell_2$ are bounded maps, and condition $(\star)$ in Theorem \ref{Thm-charact} for $T=I_{Z_2}$ implies that the map $I-\KP \circ \KP^{-1}: \ell_f^* \to \ell_2$ is also bounded.

The space $Z_2$ also appears as the derived space obtained by complex interpolation of the scale $(\ell_\infty, \ell_1)$ at $1/2$. This means that if $\mathcal C$ is the corresponding Calder\'on space, $\delta: \mathcal C\to \ell_\infty$ denotes the evaluation map at $1/2$ and $\delta': \mathcal C\to \ell_\infty$ denotes the evaluation of the derivative map at $1/2$ then $Z_2$ is isomorphically the quotient space $\mathcal C/( \ker \delta \cap \ker\delta')$, which means that (isomorphically) \begin{equation}\label{derived}Z_2 = \{(\omega, x)\in \ell_\infty\times \ell_2: \exists f\in \mathcal C: f(1/2)=x \quad\mathrm{and}\quad f'(1/2)=\omega\}\end{equation} endowed with the natural quotient norm. We will use this approach in Section \ref{opcal}.

The space $Z_2$ is not only isometric to its dual, but one singular feature of $\KP$ is its ``self-duality", in the sense that  $\KP^* = - \KP$ \cite[Theorem 5.1]{kaltpeck}, which is reflected in
$$
\left|\langle \KP x,y\rangle- \langle x,\KP y \rangle \right| \leq 2 \|x\|_2 \|y\|_2.
$$
This implies that $Z_2^*= \{(\omega^*, x^*) \in \ell_\infty\times \ell_2: \omega^* - \KP^* x^*\in \ell_2\}$  is the dual of $Z_2$ with  duality formula
\begin{equation}\label{duality}  \left \langle (\omega^*, x^*) , (\omega, x) \right \rangle = \langle \omega^*, x\rangle + \langle \omega, x^*\rangle,\end{equation}
and $D:Z_2\to Z_2^*$ given by $D(\omega, x)= (-\omega, x)$ is a bijective  isometry.

\subsection{Facts about operators on $Z_2$} The knowledge about operators on $Z_2$ is even scarcer than that about $Z_2$ itself and can be summarized in two results:

\begin{theorem}\label{Kalton} \cite[Theorems 7 and 8]{kaltsym} Let $S\in \mathfrak L(Z_2)$ and $T\in \mathfrak L(Z_2,Y)$.
\begin{itemize}
\item[(1)] If $S$ is not strictly singular then there exists a subspace $E$ of $Z_2$ isomorphic to $Z_2$ such that $S|_E$ is an isomorphism and $S(E)$ is complemented in $Z_2$ (hence $E$ also is complemented);
\item[(2)] If $T$ is not strictly singular then there exists a complemented subspace $F$ of $Z_2$ isomorphic to $Z_2$ such that $T|_F$ is an isomorphism.
\end{itemize}
\end{theorem}

\begin{prop}\label{beli} \cite[Lemma 16.15]{BeLi} A scalar $2\times 2$ matrix
$A= \left(\begin{array}{cc}
\alpha & \beta\\
\delta & \gamma \\
\end{array}
\right)$
defines an operator in $\mathfrak L(Z_2)$ in the obvious way $A(e_n, e_m) = (\alpha e_n + \beta e_m, \delta e_n + \gamma e_m)$ if and only if $\alpha=\gamma$ and $\delta=0$.\end{prop}

\section{Preliminaries} \label{sect:prelim}
\subsection{General operator theory}
An \emph{operator ideal} \cite{pietsch} is a subclass $\mathcal{A}$ of the class $\mathfrak L$ of bounded operators between Banach spaces such that finite range operators belong to $\mathcal{A}$, $\mathcal{A}+ \mathcal{A} \subset \mathcal{A}$ and {$\mathfrak L\,\mathcal{A}\,\mathfrak L\subset \mathcal{A}$.}

Let $X$ and $Y$ be Banach spaces. An operator  $T\in \mathfrak L(X,Y)$ is \emph{strictly singular} if no restriction of $T$ to an  infinite dimensional subspace is an isomorphism; $T$ is \emph{strictly cosingular} if $q_NT$ is never surjective when $q_N$ is the quotient map onto an infinite dimensional quotient $Y/N$. The classes $\mathfrak S$ of strictly singular operators and $\mathfrak C$ of strictly cosingular operators are closed operator ideals \cite[1.9 and 1.10]{pietsch}. An operator $T\in \mathfrak L(X,Y)$ is \emph{upper semi-Fredholm,} $T\in \Phi_+$, if it has closed range and finite dimensional kernel; it is \emph{lower semi-Fredholm,} $T\in \Phi_-$, if its range is finite codimensional (hence closed), and $\Phi=\Phi_+\cap\Phi_-$ is the class of Fredholm operators. Also $T$ is \emph{inessential}, denoted $T\in \mathfrak{In}$, if $I_X-AT$ is a Fredholm operator for all $A\in\mathfrak L(Y,X)$ or, equivalently, $I_Y-TA$ is Fredholm for all $A\in\mathfrak L(Y,X)$. Introduced by Kleinecke \cite{Kleinecke}, $\mathfrak {In}$ is a closed operator ideal containing both $\mathfrak S$ and $\mathfrak C$.

\subsection{Exact sequences and quasilinear maps} Let $X$ and $Y$ be quasi-Banach spaces with quasi-norms $\|\cdot\|_X$ and $\|\cdot\|_Y$.
We suppose that $Y$ is a subspace of some vector space $\Sigma$. A map $\Omega:X\to \Sigma$ is called \emph{quasilinear from $X$ to $Y$} with ambient space $\Sigma$ and denoted $\Omega: X\lop Y$ if it is homogeneous and there exists a constant $C$ so that for each $x,z\in X$,
$$
\Omega(x+z)-\Omega x-\Omega z\in Y\quad \textrm{and}\quad
\|\Omega(x+z)-\Omega x-\Omega z\|_Y \leq C(\|x\|_X +\|z\|_X).
$$

A quasilinear map $\Omega:X\lop Y$ is said to be:
\begin{itemize}
\item \emph{bounded} if there exists a constant
$D$ so that $\Omega x \in Y$ and $\|\Omega x\|_Y\leq D\|x\|_X$ for each $x\in X$.
\item \emph{trivial} if there exists a linear map $L:X\To \Sigma$ so that $\Omega - L: X\To Y$ is bounded.
\end{itemize}

A quasilinear map $\Omega: X\lop Y$ generates an exact sequence $0\to Y \to Z \to X \to 0$ (namely, a diagram formed by Banach spaces and continuous operators so that the kernel of each of them coincides with the image of the previous one), as follows: $Z =\{(\beta,x)\in \Sigma\times X : \beta-\Omega x\in Y\}
$ endowed with the quasi-norm $\|(\beta,x)\|_\Omega=\|\beta-\Omega x\|_Y+\|x\|_X$, with inclusion $y\To (y,0)$ and quotient map $(\beta,x)\To x$. The space $Z$ is sometimes called a \emph{twisted sum of $Y$ and $X$} and denoted $Y\oplus_\Omega X$. If $\Omega$ is bounded then $Y\oplus_\Omega X=Y\times X$ and $\|y-\Omega x\|_Y +\|x\|_X$
and $\|y\|_Y+\|x\|_X$ are equivalent quasi-norms on this space. If $\Omega$ is trivial then
$Y\oplus_\Omega X$ is isomorphic to $Y\oplus X$.

The general theory of twisted sums developed in \cite{kaltpeck} establishes a correspondence between exact sequences $0\To Y\To Z \To X \To 0$ of quasi Banach spaces and quasilinear maps $\Omega: X\lop Y$. The quasilinear map generating
$0\To \ell_2\To Z_2 \To \ell_2 \To 0$ is $\KP$.

\section{Matrix representation of operators on $Z_2$}\label{sect:matrix}
The space $Z_2$ admits two canonical descriptions as a twisted sum space: (\ref{ip-seq}) and (\ref{jq-seq}). Operators on $Z_2$ can be represented by a matrix $\left(\begin{array}{cc}
\alpha & \beta \\
\delta & \gamma \\
\end{array}\right),$
whose entries $\alpha, \beta, \delta$ and $\gamma$ are linear maps between sequence spaces and depend on whether one is using
(\ref{ip-seq}) or (\ref{jq-seq}). Unless specified otherwise, we will always refer to (\ref{ip-seq}); in which case the matrix above acts as
$$
\left(\begin{array}{cc}
    \alpha & \beta \\
    \delta & \gamma \\
\end{array}\right) (\omega,x)=
(\alpha\omega+ \beta x, \delta\omega +\gamma x).
$$
This same operator using (\ref{jq-seq}) has representing matrix $\left(
\begin{array}{cc}
\gamma & \delta\\
\beta & \alpha \\
\end{array}\right)$ and acts in the form
$$
\left(\begin{array}{cc}
    \gamma & \delta \\
    \beta & \alpha \\
\end{array}\right) (x, \omega)=
(\gamma x + \delta\omega , \beta x + \alpha\omega).$$

We begin with some necessary conditions:
\begin{lemma}\label{neces} \emph{[Necessary conditions]} Let $\left(
\begin{array}{cc}
\alpha & \beta \\
\delta & \gamma \\
\end{array}
\right)$ be a bounded operator $T$ on $Z_2$. The following conditions are satisfied:
\begin{itemize}
\item[(d)] $\quad\quad \delta: \ell_2\To \ell_2$ is bounded.
\item[(g)] $\quad\quad \gamma: \ell_f \To \ell_2$ and $(\beta- \KP\circ \gamma): \ell_f\To \ell_2$ are bounded.
\item[(d+gK')] $\quad\quad \delta + \gamma\circ \KP^{-1} : \ell_f^*\To \ell_2$ is a bounded map.
\item[(g+dK)] $\quad\quad \gamma + \delta\circ \KP: \ell_2\To \ell_2$ is a bounded map.
\end{itemize}\end{lemma}
\begin{proof} Recall that $i$, $p$, $j$ and $q$ are the maps appearing in the exact sequences (\ref{ip-seq}) and (\ref{jq-seq}).\medskip


(d) For $y\in\ell_2$, $pTi\,y= pT(y,0)= p(\alpha y,\delta y)= \delta y$.
\medskip

(g) For $x\in\ell_f$,
$$\|Tjx\|_{Z_2}= \|T(0,x)\|_{Z_2}= \|(\beta x,\gamma x) \|_{Z_2}= \|(\beta- \KP\circ \gamma)x\|_2+ \|\gamma x\|_2 \leq \|Tj\| \cdot\|x\|_{\ell_f}.$$
Thus $\gamma: \ell_f\To \ell_2$ and $\beta- \KP\circ \gamma: \ell_f\To \ell_2$ are bounded.\medskip

(d+gK') A bounded lifting $L_q$
for $q$ is given  by $L_q \omega =(\omega, \KP^{-1} \omega)$. Then for every $\omega \in\ell_f^*$, $\|pT L_q \omega\|_2 =\|(\delta +\gamma\circ\KP^{-1})\omega\| \leq \|pT\|\cdot\|L_q\|\cdot\|\omega\|_{\ell_f}$. Hence $\delta+ \gamma\circ \KP^{-1}: \ell_f^*\To \ell_2$ is bounded.\medskip

(g+dK) A bounded lifting $L_p$ for $p$ is given by $L_p y=(\KP y,y)$ Then for each $y\in\ell_2$, we have  $\|pT L_p y\|_2 = \|(\gamma  +\delta \circ\KP)y\| \leq \|pT\|\cdot\|L_p\|\cdot\|y\|_2$. Hence $\gamma+ \delta\circ \KP: \ell_2\To \ell_2$ is bounded.
\end{proof}

Now we characterize the bounded operators on $Z_2$ .

\begin{theorem}\label{Thm-charact}
The operator $T = \left(\begin{array}{cc}
\alpha & \beta \\
\delta & \gamma \\
\end{array}\right)$
is bounded on $Z_2$  if and only if the four necessary conditions in Lemma \ref{neces} hold as well as
$$(\star) \quad\quad \quad\alpha + \beta\circ \KP^{-1}- \KP\left(\delta+ \gamma\circ \KP^{-1}\right): \ell_f^*\To\ell_2 \quad\mathrm{is\; a\; bounded\; map}.
$$
\end{theorem}
\begin{proof}
Condition ($\star$) is necessary: if $T$ is bounded then $\|\alpha \omega + \beta x - \KP(\delta \omega+ \gamma x)\|_2\leq\|T\| \|(\omega, x)\|_{Z_2}$ and the choice
$(\omega, \KP^{-1}\omega)= L_q\omega$ yields
$$
\|\alpha \omega + \beta\circ \KP^{-1} \omega- \KP(\delta\omega+ \gamma\circ \KP^{-1}\omega)\|_2\leq \|L_q\|\|T\|  \|\omega\|_{\ell_f^*}.
$$

Conversely, we will show that there exists a constant $C>0$ so that for each $(\omega,z) \in Z_2$ we have $(\alpha\omega+ \beta x, \delta \omega+ \gamma x)\in Z_2$ and $\|(\alpha\omega+ \beta x,\delta\omega+ \gamma x)\|_{Z_2} \leq C\|(\omega,x)\|_{Z_2}$. We need to show that

(1) $\delta\omega+ \gamma x\in \ell_2$,

(2) $\alpha\omega+ \beta x-\KP(\delta \omega+ \gamma x)\in\ell_2$\quad (hence $\alpha\omega+ \beta x\in\ell_\infty$), and

(3) $\|\delta \omega+\gamma x\|_2+ \|\alpha \omega + \beta x - \KP \left( \delta \omega + \gamma x \right)\|_2\leq C \|(\omega, x)\|_{Z_2}.$
\medskip

Recall that there exist $M>0$ such that $\|x- \KP^{-1}\omega \|_{\ell_f}+ \|\omega\|_{\ell_f^*}\leq M \|(\omega,x) \|_{Z_2}$ for each $(\omega, x) \in Z_2$.

(1) Observe that $\omega -\KP x\in\ell_2$ and  $x-\KP^{-1}\omega \in\ell_f$, and by Lemma \ref{neces}, the maps $\gamma+ \delta\circ\KP: \ell_2\To \ell_2$ and $\delta +\gamma\circ\KP^{-1}: \ell_f^*\To \ell_2$ are bounded. Thus
$$
\delta\omega+ \gamma x= \frac{1}{2} \Big( \delta(\omega-\KP x)+ \gamma(x-\KP^{-1}\omega) + (\gamma+\delta \circ\KP)x+ (\delta+\gamma\circ \KP^{-1})\omega\Big)\in \ell_2
$$
and
\begin{eqnarray*}
\|\delta\omega+ \gamma x\|_2 &\leq& \frac{1}{2} \Big(\|\delta\|\|\omega-\KP x\|_2+ \|\gamma\| \|x-\KP^{-1}\omega\|_{\ell_f} + \|\gamma+\delta \circ\KP\|\|x\|_2\\
&+& \|\delta+ \gamma\circ \KP^{-1}\| \|\omega\|_{\ell_f^*}\Big)\\
&\leq& \frac{1}{2} \Big(\|\delta\|+ \|\gamma\|M+ \|\gamma+\delta \circ\KP\|+ \|\delta+\gamma\circ \KP^{-1}\|M\Big) \|(\omega, x)\|_{Z_2}.
\end{eqnarray*}

To prove (2) we decompose $\alpha \omega + \beta x - \KP \left( \delta \omega + \gamma x \right)$ in three pieces:
\begin{eqnarray*}
&& \alpha \omega + \beta\circ \KP^{-1} \omega-\KP \left(\delta+ \gamma\circ \KP^{-1}\right)\omega\\
&+& \beta\left(x- \KP^{-1} \omega \right)- \KP\left(\gamma x- \gamma \circ \KP^{-1} \omega \right)\\
&+& \KP\left(\gamma x - \gamma\circ \KP^{-1} \omega \right)
+ \KP\left( \delta + \gamma\circ \KP^{-1}\right) \omega - \KP\left( \delta \omega + \gamma x\right).
\end{eqnarray*}

The first piece is bounded by $(\star)$; and for the third piece, note that  $\KP: \ell_2 \curvearrowright \ell_2$ is quasilinear, hence $\KP(x+y)-\KP x-\KP y\in\ell_2$ and
$$\|\KP(x+y)-\KP x-\KP y\|_2 \leq \|\KP\| (\|x\|_2+\|y\|_2)$$
for $x,y\in\ell_2$. Moreover  $\gamma$ and $\delta+\gamma\circ \KP^{-1}$ are both bounded from $\ell_f$ into $\ell_2$. Thus
\begin{eqnarray*}
&&\|\KP(\gamma x- \gamma\circ \KP^{-1} \omega)+ \KP(\delta\omega + \gamma\circ \KP^{-1} \omega) - \KP(\delta \omega + \gamma x)\|_2\\
&\leq& \|\gamma x - \gamma\circ \KP^{-1} \omega \|_2 + \|\delta \omega + \gamma\circ \KP^{-1} \omega\|_2\\
&\leq& \|\gamma\|\| x - \KP^{-1} \omega \|_{\ell_f} + \|\delta + \gamma\circ \KP^{-1}\| \|\omega\|_{\ell_f^*}\\
&\leq& (\|\gamma\|+ \|\delta + \gamma\circ \KP^{-1}\|) M\|(\omega,x)\|_{Z_2}.
\end{eqnarray*}

For the second piece, since $x - \KP^{-1} \omega \in \ell_f$, $(\beta - \KP\circ\gamma):\ell_f\To \ell_2$ is bounded and
$\beta\left(x - \KP^{-1} \omega\right) - \KP\left(\gamma x - \gamma\circ \KP^{-1} \omega \right)$ is equal to
$$\beta\left( x - \KP^{-1} \omega\right) - \KP\circ\gamma \left(x - \KP^{-1} \omega \right)\\
= (\beta - \KP\circ\gamma) \left(x - \KP^{-1} \omega \right),$$
one has
\begin{eqnarray*}
\|(\beta - \KP\circ\gamma) \left(x - \KP^{-1} \omega \right)\|_2
&\leq&  \|(\beta - \KP\circ\gamma) \| \|x-\KP^{-1} \omega \|_{\ell_f}\\ &\leq& M\|(\beta - \KP\circ\gamma)\| \|(\omega, x)\|_{Z_2}.
\end{eqnarray*}

(3) clearly follows from the arguments in the proof of (1) and (2). \end{proof}

Condition ($\star$) for $T=I_{Z_2}$ implies that $I- \KP\circ\KP^{-1}: \ell_f^*\To\ell_2$ is bounded, from which we can derive the Benyamini-Lindenstrauss characterization of Proposition \ref{beli}: if $\left(\begin{array}{cc}
\alpha & \beta \\
\delta & \gamma \\
\end{array}\right)$
is a bounded operator on $Z_2$ with $\alpha, \beta, \gamma,\delta$ scalars then the boundedness of $\gamma + \delta \KP: \ell_2\to \ell_2$ implies that $\delta=0$ while the boundedness of $\alpha + \beta \KP^{-1} - \gamma \KP \circ \KP^{-1}: \ell_f^* \to \ell_2$ yields that, since $\beta \KP^{-1}$ is also bounded for any scalar $\beta$ and $ \alpha - \gamma + \gamma (I - \KP\circ \KP ^{-1})$ is bounded then also $\alpha - \gamma: \ell_f^*\to \ell_2$ is bounded, and thus $\alpha = \gamma$.

Apart from ($\star$) and the necessary conditions in Lemma \ref{neces}, we have  some additional ones that were not needed in the proof of Theorem \ref{Thm-charact}:

\begin{lemma}\label{addit-Thm}
If  $T = \left(
\begin{array}{cc}
\alpha & \beta \\
\delta & \gamma \\
\end{array}
\right)\in \mathfrak L(Z_2)$  then the following conditions are satisfied:
\begin{itemize}
\item[(a)] $\quad\quad \alpha:\ell_2\To \ell_f^*$ is a bounded operator.
\item[(b)] $\quad\quad \beta: \ell_f \To \ell_f^*$ is a bounded operator.
\item[(c$_0$)]\quad  $(\alpha- \KP\circ \delta): \ell_2\To \ell_2$ is bounded.
\item[(c$_1$)]\quad $\alpha\circ\KP+ \beta: \ell_2\To\ell_f^*$ is bounded.
\item[(c$_2$)]\quad $\alpha+\beta\circ \KP^{-1}: \ell_f^*\To \ell_f^*$ is bounded.
\item[(c$_3$)]\quad $\alpha\circ\KP+\beta -\KP(\delta \circ\KP +\gamma):\ell_2\To \ell_2$ is bounded.
\item[(c$_4$)]\quad $\gamma  - \KP^{-1} \circ \beta: \ell_f \to \ell_f$ is bounded.
 \end{itemize}
\end{lemma}
\begin{proof}
(a) and (b) follow from $\alpha=qTi$ and $\beta=qTj$.

(c$_0$) For $y\in\ell_2$,
$\|Tiy\|_{Z_2}= \|T(y,0)\|_{Z_2}= \|(\alpha y, \delta y) \|_{Z_2}= \|(\alpha- \KP\circ \delta)y\|_2+ \|\delta y\|_2\leq \|Ti\|\cdot\|y\|_2$. Therefore $(\alpha- \KP\circ \delta): \ell_2\To \ell_2$ is bounded.

(c$_1$) For each $y\in\ell_2$, $\|qTL_p y\|=
\|qT(\KP y,y)\|_{Z_2}= \|(\alpha\circ\KP+ \beta)y\|_{\ell_f^*}\leq \|T\|\|L_p\|\|y\|_2$.

(c$_2$) and (c$_3$) are proved in a similar way, using that $qTL_q \omega = (\alpha+ \beta\circ \KP^{-1}) \omega$ for each $\omega\in \ell_f^*$ and $\|TL_p y\|_{Z_2}=  \|(\alpha \circ \KP y,\delta\circ\KP y) \|_{Z_2} \leq \|T\|\|L_p\|\|y\|_2$.

(c$_4$) follows from the continuity of $Tj$. \end{proof}

Let us see some applications.

\begin{prop}\label{additional} Let $T = \left(
\begin{array}{cc}
\alpha & \beta \\
\delta & \gamma \\
\end{array}
\right)$ be an operator on $Z_2$ .
\begin{itemize}
\item[(d*)] If $\gamma: \ell_2\to \ell_2$ is bounded then $\delta: \ell_f^*\To \ell_2$ and $\delta\circ \KP: \ell_2\To \ell_2$ are also bounded.
\item[(d')] If $\alpha: \ell_2\to \ell_2$ is bounded then $\delta: \ell_2\To \ell_f$ is bounded. Hence  $\KP\circ \delta: \ell_2\To \ell_2$ is bounded.
\item[(d**)] If $\gamma: \ell_f\to \ell_f$ is bounded then $\KP^{-1}\circ \beta: \ell_f\To \ell_f$ is also bounded.
\item[(d'')] If $\alpha: \ell_f^*\to \ell_f^*$ is bounded then $\beta \circ \KP^{-1}: \ell_f^*\To \ell_f^*$ is bounded.
\end{itemize}
\end{prop}

\begin{proof}
(d*) Since $\KP^{-1}:\ell_f^*\To \ell_2$ and $\gamma: \ell_2\To \ell_2$ are bounded, so is $\gamma\circ  \KP^{-1}: \ell_f^*\To \ell_2$. By (d+gK') in  Lemma \ref{neces}, $\delta: \ell_f^*\To \ell_2$ is bounded, hence $\delta \circ \KP: \ell_2\To \ell_2$ is also bounded.

(d') By (c$_0$) in Lemma \ref{addit-Thm}, $\alpha- \KP\circ \delta: \ell_2\To \ell_2$ is bounded. Then $\alpha: \ell_2 \To \ell_2$  bounded implies  $\KP\circ \delta: \ell_2\To \ell_2$ bounded; hence $\KP^{-1}\circ \KP\circ\delta: \ell_2\To \ell_f$ bounded; equivalently, 
$\delta: \ell_2\To \ell_f$ is bounded. For the last equivalence observe that, by the definition of domain of $\KP$ and $\KP^{-1}$, for $x\in\ell_2$, $\KP^{-1}\circ \KP x\in \ell_f \Rightarrow \KP x\in \ell_2 \Rightarrow x\in\ell_f$. The assertions (d**) and (d'') can be obtained analogously. \end{proof}

Proposition \ref{additional} yields the boundedness of $\KP\circ\delta, \delta\circ\KP, \beta \circ \KP^{-1}$ or $ \KP^{-1}\circ \beta$ but only under additional conditions that, in general, are not guaranteed. It is for that reasons surprising that one has:

\begin{prop}\label{uselemma3} Let $T = \left(
\begin{array}{cc}
\alpha & \beta \\
\delta & \gamma \\
\end{array}
\right) \in \mathfrak L(Z_2)$. Then:
 \begin{enumerate}
 \item  $\KP\circ \delta:\ell_2 \curvearrowright \ell_2$ and $\delta\circ \KP:\ell_2\curvearrowright \ell_2$ are trivial quasilinear maps.
 \item  $\KP^{-1}\circ \beta:\ell_f \curvearrowright \ell_f$ and $\beta\circ \KP^{-1}:\ell_f^*\curvearrowright \ell_f^*$ are trivial quasilinear maps.\end{enumerate} \end{prop}
\begin{proof} We prove assertion (1). Consider the pull-back  diagram
$$\xymatrix{0\ar[r]&\ell_2 \ar[r]^i &Z_2  \ar[r]^p& \ell_2\ar[r]&0\\
0\ar[r]&\ell_2 \ar@{=}[u] \ar[r] & \PB \ar[r]\ar[u]& \ell_2\ar[r]\ar[u]_\delta&0}$$

Then $\KP\circ\delta$ is a  quasilinear map generating the lower exact sequence and $\KP\circ\delta$ is trivial if and only if $\delta$ admits a bounded linear lifting $\ell_2\To Z_2$ \cite[Lemma 2.8.3]{hmbst}. Since $\alpha-\KP\circ\delta:\ell_2\to \ell_2$ and $\delta:\ell_2\To \ell_2$ are  bounded, $\hat\delta x=(\alpha x,\delta x)$ provides the lifting. Indeed, $\|(\alpha x,\delta x)\|_{Z_2}= \|(\alpha- \KP\circ\delta)x\|_2+ \|\delta x\|_2$.

Analogously, if we consider the push-out diagram
$$\xymatrix{0\ar[r]&\ell_2 \ar[d]_\delta \ar[r]^i &Z_2  \ar[d]\ar[r]^p& \ell_2\ar[r]&0\\
0\ar[r]&\ell_2 \ar[r] & \PO \ar[r]& \ell_2\ar[r]\ar@{=}[u] &0}$$
then $\delta\circ\KP$ is a quasilinear map generating the lower exact sequence which is trivial if and only if $\delta$ admits a bounded linear extension $\overline\delta:Z_2\To \ell_2$ \cite[Lemma 2.6.3]{hmbst}. In the proof of Theorem \ref{Thm-charact} it is showed that $\overline\delta(\omega,x)=\delta \omega+\gamma x$ provides such an extension.\medskip

The proof for assertion (2) follows in a similar way: to show that $\KP^{-1}\circ \beta: \ell_f \lop \ell_f$ just consider the pull-back diagram
$$\xymatrix{0\ar[r]&\ell_f \ar[r]^j &Z_2  \ar[r]^q& \ell_f^*\ar[r]&0\\
0\ar[r]&\ell_f \ar@{=}[u] \ar[r] & \PB \ar[r]\ar[u]& \ell_f\ar[r]\ar[u]_\beta&0}$$
(where $\beta$ is continuous by Proposition \ref{addit-Thm} (b)) and observe that $\omega \to (\beta \omega, \gamma \omega)$ is a bounded lifting for $\beta$ by $(c_4)$. To show that to show that $\beta\circ \KP^{-1}: \ell_f^* \lop \ell_f^*$ just consider the pushout diagram
$$\xymatrix{0\ar[r]&\ell_f \ar[r]^j \ar[d]_\beta &Z_2 \ar[d] \ar[r]^q& \ell_f^*\ar[r]&0\\
0\ar[r]&\ell_f^* \ar[r] & \PO \ar[r]& \ell_f^*\ar[r]\ar@{=}[u]&0}$$
and use $(c_2)$. \end{proof}


\section{Triangular operators} \label{sect:triangular}

An operator $T\in\mathfrak L(Z_2)$ is said to be \emph{compatible with the presentation (\ref{ip-seq})} if it satisfies $T[i[\ell_2]]\subset i[\ell_2]$. This occurs if and only if its corresponding matrix
$\left(\begin{array}{cc}
\alpha & \beta \\
0 & \gamma \\
\end{array}\right)$
is \emph{upper triangular}; namely, $\delta=0$. In that case, the diagram
$$
\xymatrix{0\ar[r]&\ell_2\ar[d]_\alpha \ar[r]^i &Z_2  \ar[d]^T\ar[r]^p& \ell_2\ar[r]\ar[d]^\gamma&0\\
0\ar[r]&\ell_2 \ar[r]^i &Z_2 \ar[r]^p& \ell_2\ar[r]&0}
$$
is commutative. The operator
$T=\left(\begin{array}{cc}
\alpha & \beta \\
\delta & \gamma \\
\end{array}\right)$
above is said to be \emph{compatible with the presentation (\ref{jq-seq})} if it satisfies $T[j[\ell_f]]\subset j[\ell_f]$, and this occurs if and only if its matrix is \emph{lower triangular}; namely, $\beta=0$. In that case, the diagram
$$\xymatrix{0\ar[r]&\ell_f\ar[d]_\gamma \ar[r]^j  &Z_2  \ar[d]^T\ar[r]^q& \ell_f^* \ar[r]\ar[d]^\alpha&0\\
0\ar[r]&\ell_f \ar[r]^j &Z_2 \ar[r]^q& \ell_f^*\ar[r]&0.}
$$is commutative. Consequently, $T:Z_2\to Z_2$ is compatible with both presentations if and only if it is \emph{diagonal}, in the sense that its matrix is $\left(
\begin{array}{cc}
\alpha & 0\\
0 & \gamma \\
\end{array}\right).$ We show next that the additional conditions of Proposition \ref{additional} are always satisfied for $T$ upper triangular.

\begin{lemma}\label{ncfutm}\emph{[Necessary conditions]} If $T = \left(
\begin{array}{cc}
\alpha & \beta \\
0 & \gamma \\
\end{array}
\right): Z_2\To Z_2$ is a bounded operator then
\begin{itemize}
\item[(a')] $\quad\quad\alpha: \ell_2\To \ell_2$ is a bounded operator.
\item[(g')] $\quad\quad\gamma: \ell_2 \To \ell_2$ is a bounded operator.
\item[(k)] $\quad\quad\alpha - \gamma$ is compact.
\end{itemize}\end{lemma}
\begin{proof}
Since $\delta=0$, $\alpha: \ell_2\To \ell_2$ and $\gamma: \ell_2 \To \ell_2$ are bounded by (c$_0$) in Lemma \ref{addit-Thm} and (g+dK) in Lemma \ref{neces}.

(k) was proved in \cite[Corollary 5.9]{ccfm}.
%
\end{proof}

The necessary conditions become sufficient if we add (see Theorem 1) that
$$\alpha + \beta\circ \KP^{-1}- \KP\circ  \gamma\circ \KP^{-1}: \ell_f^* \To \ell_2\quad \textrm{is a bounded map.}
$$

Equivalently, $\alpha \omega +\beta x-\KP\circ\gamma x$ is a bounded map $Z_2\To \ell_2$ (recall that $(\omega, x)\in Z_2$ means that $\omega - \KP x\in \ell_2$ or $x - \KP^{-1} \omega \in \ell_f$). Therefore

\begin{theorem}
An operator
$T = \left(\begin{array}{cc}
\alpha & \beta \\
0 & \gamma \\
\end{array} \right)$
is bounded on $Z_2$ if and only if $\alpha, \gamma\in \mathfrak L(\ell_2)$ and the maps $\beta- \KP\circ \gamma: \ell_f \To \ell_2$ and  $\alpha+ \beta \circ\KP^{-1}- \KP\circ \gamma\circ \KP^{-1}:\ell_f^* \To\ell_2$ are bounded.
\end{theorem}

A few variations of the previous result are possible:
\begin{prop}\label{cases}$\;$
\begin{itemize}
\item[(a)] $S=\left(\begin{array}{cc}
\alpha & 0 \\
\delta & 0 \\
\end{array}\right) \in \mathfrak L(Z_2)$
if and only if $\KP\circ \delta-\alpha: \ell_f^* \To \ell_2$ is bounded. If so, $S\in \mathfrak S$. In particular, $\left(\begin{array}{cc}
\alpha & 0 \\
0 & 0 \\
\end{array}\right)\in \mathfrak L(Z_2)$
if and only if $\alpha \in \mathfrak L(\ell_f^*,\ell_2)$; and
$\left(\begin{array}{cc}
0 & 0 \\
\delta & 0 \\
\end{array}
\right)\in \mathfrak L(Z_2)$
if and only if $\delta \in \mathfrak L(\ell_f^*,\ell_f)$.

\item[(b)] $R = \left(\begin{array}{cc}
0 & \beta \\
0 & \gamma \\
\end{array}
\right)\in \mathfrak L(Z_2)$
if and only if $\KP\circ \gamma$ is trivial and  $\KP\circ \gamma- \beta: \ell_2\To \ell_2$ is bounded. If so, $R\in \mathfrak S$. In particular, $\left(\begin{array}{cc}
0 & \beta \\
0 & 0 \\
\end{array}
\right)\in \mathfrak L(Z_2)$
if and only if $\beta \in \mathfrak L(\ell_2,\ell_2)$; and
$\left(\begin{array}{cc}
0 & 0 \\
0 & \gamma \\
\end{array}\right)\in \mathfrak L(Z_2)$
if and only if $\gamma \in \mathfrak L(\ell_2,\ell_f)$.
\item[(c)] $T=\left(\begin{array}{cc}
\alpha & \beta \\
0 & 0 \\
\end{array}\right) \in \mathfrak L(Z_2)$ if and only if $\alpha \circ \KP + \beta: \ell_2\To \ell_2$ is bounded. If so, $T\in\mathfrak S$.
\medskip

\item[(d)] $M=\left(\begin{array}{cc}
0 & 0 \\
\delta & \gamma \\
\end{array}\right) \in \mathfrak L(Z_2)$
if and only if $\delta \circ \KP + \gamma : \ell_2\To \ell_f$ is bounded. If so, $M\in\mathfrak S$.
\end{itemize}
\end{prop}

\begin{proof}
Use Proposition \ref{uselemma3}. (a) Suppose $S\in \mathfrak L(Z_2)$. Since $S(\omega,x)= (\alpha\omega, \delta \omega)$, $S$ factors through $\ell_f^*$; in fact $S=Bq$ with $B:\ell_f^*\to Z_2$ given by $B\omega = (\alpha\omega, \delta \omega)$, and $B$ is bounded because $\KP\circ \delta$ is trivial with $\KP\circ \delta- \alpha$ bounded; hence $\omega\To (\alpha\omega, \delta\omega)$ is a continuous lifting $\ell_f^* \To Z_2$ for $\delta: \ell_f^* \To \ell_2$ as in the diagram
$$\xymatrix{
0\ar[r]& \ell_2\ar[r]& Z_2\ar[r]^p &\ell_2 \ar[r]& 0\\
&&&\ell_f^* \ar[ul]^{(\alpha, \delta)} \ar[u]_\delta&}$$
Moreover, since $p$ is strictly singular, so is $S$. Regarding the ``in particular" assertion, it is clear that if $\delta \in \mathfrak L(\ell_f^*, \ell_f)$ then $\KP\circ \delta: \ell_f^*\to \ell_2$ is bounded and, by (a), the corresponding $S$ is bounded. Conversely, if $S$ is bounded on $Z_2$ then $\delta[\ell_f^*]\subset \ell_f$. Moreover, we have $\|\delta \omega\|_{\ell_f} = \|(0, \delta \omega) \|_{Z_2} \leq \|S\| \|(\omega, x) \|_{Z_2}$ and therefore $\delta: \ell_f^* \to \ell_f$ is bounded. The proofs of (b,c,d) are analogous. \end{proof}

Johnson, Lindenstrauss and Schechtman \cite{jls} asked whether every operator $T: Z_2 \to Z_2$ is a strictly singular perturbation of an operator sending $i[\ell_2]$ into itself.
We could also consider operators sending $j[\ell_f]$ into itself and formulate the corresponding conjecture.
One has the following partial result:
\begin{prop}
Let $T=\left(\begin{array}{cc}
\alpha & \beta \\
\delta & \gamma \\
\end{array}\right) \in \mathfrak L(Z_2)$
 \begin{itemize}
 \item[(1)] If $\delta \in \mathfrak L(\ell_f^*, \ell_f)$ then $T$ is a strictly singular perturbation of an upper triangular operator.
 \item[(2)] If $\beta \in \mathfrak L(\ell_2, \ell_2)$ then $T$ is a strictly singular perturbation of a lower triangular operator.\end{itemize}
\end{prop}
\begin{proof}
(1) follows from
$\left(\begin{array}{cc}
\alpha & \beta \\
\delta & \gamma \\
\end{array}\right) =
\left(\begin{array}{cc}
\alpha & \beta \\
0 & \gamma \\
\end{array}\right) +
\left(\begin{array}{cc}
0 &  0 \\
\delta & 0 \\
\end{array}
\right)$
and Proposition \ref{cases} (a), and (2) follows from
$\left(\begin{array}{cc}
\alpha & \beta \\
\delta & \gamma \\
\end{array}\right) =
\left(\begin{array}{cc}
\alpha & 0 \\
\delta & \gamma \\
\end{array}\right) +
\left(\begin{array}{ccc}
0 &  \beta \\
0 & 0 \\
\end{array}
\right)$ and Proposition \ref{cases}(b). \end{proof}

\section{Complemented copies of $Z_2$ in $Z_2$} \label{sect:involution}

Kalton \cite{kaltsym} considered the continuous alternating bilinear form
$\Omega: Z_2\times Z_2\to\R$ given by
$$
\Omega \left((\omega_1, x_1), (\omega_2,x_2) \right\rangle = \langle\omega_1, x_2\rangle- \langle \omega_2, x_1\rangle,
$$
and for every $T\in \mathfrak{L}(Z_2)$ he defined an operator $T^+:Z_2\To Z_2$ by $\Omega(T^+y,z)= \Omega(y,Tz)$ for $y,z\in Z_2$. Note that $(Dx)y= \Omega(x,y)$ for $y,z\in Z_2$. Hence $T^+=D^{-1}T^* D\in \mathfrak{L}(Z_2)$, where $T^*\in \mathfrak{L}(Z_2^*)$ is the conjugate operator of $T$.
Indeed, for $x,y\in Z_2$ we have
$$(DT^+x)y= \Omega(T^+x,y)=  \Omega(x,Ty)= Dx(Ty)= (T^*Dx)y;$$ hence $DT^+=T^*D$.
Moreover, the map $T\To T^+$ is an \emph{involution} on $\mathfrak{L}(Z_2)$ \cite[Definition 11.14]{Rudin}.\medskip

Since $D$ is a bijective isometry, most of the properties of $T^+$ coincide with those of $T^*$. Namely, $\|T\|=\|T^+\|$, $R(T)$ is closed if and only if $R(T^+)$ is so, $T$ is an isomorphism into if and only if $T^+$ is surjective, and $T\in\Phi_-$ if and only if $T^+\in\Phi_+$. In particular, $T^+=T$ and $T\in\Phi_+$ imply $T\in\Phi$.
\medskip

It is easy to check that if\;
$T=\left(\begin{array}{cc}
\alpha & \beta \\
\delta & \gamma \\
\end{array}
\right)$\; then\;
$T^*= \left(\begin{array}{cc}
\gamma^* & \beta^* \\
\delta^* & \alpha^* \\
\end{array}\right)$\;
and\;
$T^+= \left(\begin{array}{cc}
\gamma^* & -\beta^* \\
-\delta^* & \alpha^* \\
\end{array}\right).$

\adef \emph{\cite{kaltsym}}
A subspace $E$ of $Z_2$ is said to be \emph{isotropic} when $\Omega(u,v)= 0$ for every $u,v\in E$. An operator $T\in \mathfrak{L}(Z_2)$ is \emph{isotropic} when its range is isotropic; equivalently, when $T^+T=0$.\zdef

Clearly $i[\ell_2]$, $j[\ell_f]$ and $\{(x,x)\in Z_2 : x\in \ell_f\}$ are isotropic subspaces of $Z_2$. Moreover, the operator $i\, p=\begin{pmatrix} 0 & I\\ 0 & 0 \end{pmatrix}\in \mathfrak{L}(Z_2)$ is isotropic, non-compact, strictly singular and  strictly cosingular, with $R(i\, p)=N(i\, p)$ and $(i\, p)^+=-i\, p$.

\begin{lemma}
For each $T\in \mathfrak{L}(Z_2)$, $R(T)\cap N(T^+)$ is an isotropic subspace.
\end{lemma}
\begin{proof}
Let $(\omega_1,x_1), (\omega_2,x_2) \in R(T)\cap N(T^+)$. Then $(\omega_2,x_2)= T(\omega,x)$ for some $(\omega,x)\in Z_2$ and
$$\Omega \left((\omega_1, x_1), (\omega_2,x_2) \right\rangle = \Omega \left((\omega_1, x_1), T(\omega,x) \right\rangle = \Omega \left(T^+ (\omega_1, x_1), (\omega,x) \right\rangle = 0$$
concluding the proof.
\end{proof}

\begin{prop}\label{restrict}
Let $T\in \mathfrak{L}(Z_2)$.
\begin{itemize}
\item[(a)] If $T^+T$ is strictly singular then so is $T$.
\item[(b)] If $Ti$ is strictly singular then so is $T$.
\item[(c)] If $Ti\in\Phi_+$ then $T\in\Phi_+$.
\end{itemize}
\end{prop}
\begin{proof}
(a) and (b) are \cite[Theorem 9 and Lemma 5]{kaltsym}.

(c) Suppose that $T\not\in\Phi_+$. Then there exists an infinite dimensional subspace $M \subset Z_2$ such that $T|_M$ is compact. Since $p$ is strictly singular, we can find a normalized basic sequence $(x_n)$ in $M$ and  sequences $(x^*_n)$ in $Z_2^*$ with $C=\sup_n \|x^*_n\|<\infty$ and $(y_n)$ in $i[\ell_2]=N(p)$ with $\|x_n-y_n\|< 2^{-n}/C$. Then the expression
$Kx= \sum_{i=1}^\infty x^*_i(x) (x_i-y_i)$ defines a compact operator $K\in \mathfrak{L}(Z_2)$ with $\|K\|<1$ such that, denoting by $N$ the closed subspace generated by $(x_n)$, we have $(I-K)[N]\subset i[\ell_2]$. We claim that $Ti|_N$ is compact, hence $Ti\not\in\Phi_+$.
Indeed, if $(z_n)$ is a bounded sequence in $N$, then $(T(I-K)z_n)= (Tz_n-TKz_n)$ has a convergent subsequence. Since $TK$ is compact, $(Tz_n)$ has a convergent subsequence.
\end{proof}

\begin{lemma}
Let $T\in \mathfrak{L}(Z_2)$. If $T\in\Phi_+$  then $T^+T\in\Phi$.
\end{lemma}
\begin{proof}
Since $(T^+T)^+=T^+T$, it is enough to show that $T^+T\in\Phi_+$. Suppose that $T^+T\not\in\Phi_+$. Then $T^+Ti\not\in\Phi_+$; hence there exists a normalized block basis sequence $(w^{(n)})$ in $\ell_2$ such that, if we denote by $i_w: \ell_2\to \ell_2$ the isometric embedding defined by $i_w e_n= w^{(n)}$, then $T^+Ti i_w$ is compact. Let $W$ be the block operator associated to the sequence $(w^{(n)})$ in \cite[Section 4]{kaltsym}. Since $T^+TWie_n= T^+Ti i_we_n$ for each $n\in\N$, $T^+TWi$ is compact; hence $T^+TW \in \mathfrak S$ by (b) in Proposition \ref{restrict}.
Thus $W^+T^+TW\in\mathfrak S$, hence $TW\in \mathfrak S$ by (a) in Proposition \ref{restrict}, implying $T\notin\Phi_+$.
\end{proof}

Observe that $T^+T\in\Phi_+$ implies $T\in\Phi_+$.\medskip

It is known \cite[Theorem 16.16]{BeLi} that every infinite dimensional complemented subspace of $Z_2$ contains a further complemented subspace isomorphic to $Z_2$. It is not known whether every infinite dimensional complemented subspace of $Z_2$ is isomorphic to $Z_2$ (which in particular would imply that $Z_2$ is isomorphic to its hyperplanes). We add now a new piece of knowledge:

\begin{theorem}\label{Z2inZ2}
Every subspace of $Z_2$ isomorphic to $Z_2$ is complemented.
\end{theorem}
\begin{proof}
Let $T\in \mathfrak{L}(Z_2)$ an isomorphism into with range $R(T)$. Then $T^+T\in\Phi$, hence here exists a finite codimensional subspace $N$ of $R(T)$ such that $T^+|_N$ is an isomorphism and $T^+[N]$ is finite codimensional, hence complemented; thus $N$ is complemented and so is $R(T)$.
\end{proof}

An extension of Theorem \ref{Z2inZ2} in operator terms is available now:

\begin{theorem}
Every semi-Fredholm operator on $Z_2$ has complemented kernel and range.
\end{theorem}
\begin{proof} Let $T$ be an operator on $Z_2$. If $T\in \Phi_+$ then the kernel is finite dimensional, and we can prove that $R(T)$ is complemented with the proof of Theorem \ref{Z2inZ2}. If $T\in \Phi_-$ then $R(T)$ is closed  finite codimensional and $T^*\in \Phi_+$ (in $Z_2^*$). Since $Z_2^*\simeq Z_2$, $R(T^*)$ is complemented by the first part, hence $N(T)= ^\perp R(T^*)$ is also complemented. \end{proof}

Recall from \cite{2132} that a Banach space $X$ is said to be \emph{$Y$-automorphic} if every isomorphism between two infinite codimensional subspaces of $X$ isomorphic to $Y$ can be extended to an automorphism of $Z_2$. It is clear that $\ell_2$ is $\ell_2$-automorphic and that $Z_2$ is not $\ell_2$-automorphic: indeed, since $Z_2\simeq Z_2 \oplus Z_2$, an isomorphism between the subspaces $\ell_2\oplus 0$ and $\ell_2\oplus \ell_2$ cannot be extended to an automorphism of $Z_2$.
Surprisingly,  one  has:

\begin{prop} $Z_2$ is $Z_2$-automorphic. \end{prop}
\begin{proof} As Kalton remarks in \cite[p. 110]{kaltsym}, Pe\l czy\'nski's decomposition argument shows that if $E$ is a complemented subspace of $Z_2$ and $E\oplus E\simeq E$ then $E$ is isomorphic to $Z_2$.
Suppose that $Z_2\simeq Z_2\oplus F$. By (2) of Theorem \ref{Kalton} one has $F\simeq Z_2 \oplus N$, and thus $F\oplus F\simeq F\oplus Z_2\oplus N\simeq Z_2\oplus N\simeq F$. Hence $F\simeq Z_2$.
\end{proof}

Now if $E$ is an infinite dimensional complemented subspace of $Z_2$ then
$$E\simeq Z_2\oplus E' \simeq Z_2\oplus Z_2\oplus E' \simeq Z_2 \oplus E$$
and thus $E$ is isomorphic to its $2$-codimensional subspaces since
$$E \oplus \mathbb K^2 \simeq Z_2\oplus E\oplus \mathbb K^2 \simeq Z_2 \oplus E \simeq E.$$
We conjecture that $E\oplus E\simeq Z_2$.
\medskip


Since $Z_2^*\simeq Z_2$ and $i$ is strictly cosingular \cite{kaltpeck} (hence $i^*$ strictly singular), one can easily derive the following results, by duality, from the previous ones:

\begin{prop}\label{restrict-dual}
Let $T\in \mathfrak{L}(Z_2)$.
\begin{itemize}
\item[(a)] If $T^+T$ is strictly cosingular then so is $T^+$.
\item[(b)] If $pT$ is strictly  cosingular then so is $T$.
\item[(c)] If $pT\in\Phi_-$ then $T\in\Phi_-$.
\item[(d)]  If $T^+\in\Phi_-$  then $T^+T\in\Phi$.
\end{itemize}
\end{prop}

\section{Examples of operators on $Z_2$} \label{sect:examples}

\subsection{Rank-one operators}
Given $(x^*,\omega^*)\in Z_2^*$ and $(u, v)\in Z_2$, the rank-one operator $(x^*,\omega^*)\otimes (u, v)$ acts on  $Z_2$ as follows. For each $(\omega,x)\in Z_2$,
$$
[(x^*,\omega^*)\otimes(u, v)](\omega,x)= \langle (x^*,\omega^*),(\omega,x)\rangle\cdot (u,v) = (\omega \omega^*+ x^* x)\cdot (u,v).
$$
Thus the matrix associated to $(x^*,\omega^*)\otimes (u,v)$ is
$\left(\begin{array}{cc}
\omega^*\otimes u & x^*\otimes u \\
\omega^*\otimes v & x^*\otimes v\\
\end{array}\right).$

\subsection{Nuclear operators} Fix $u^*\in Z_2^*$ and $v\in Z_2$ and obtain the matrix representation for the one-dimensional map
$u^*\otimes v$.
If $u^* =\sum a_nu_n^*$ and $v =\sum b_n v_n$, set $\mathbf u(2n-1)^* =\sum a_{2n-1} u_{2n-1}$ the ``odd" part of $u^*$
and $\mathbf u(2n)^* =\sum a_{2n} u_{2n}$ the ``even" part of $u^*$. Define in the same manner the odd and even parts $\mathbf v(2n-1)$ and $\mathbf v(2n)$ of $v$ to obtain

$$u^* \otimes v = \left(\begin{array}{cc}
\mathbf u(2n-1)^* \otimes \mathbf v(2n-1)& \mathbf u(2n)^* \otimes \mathbf v(2n-1) \\
\mathbf u(2n-1)^*\otimes \mathbf v(2n) & \mathbf u(2n)^* \otimes \mathbf v(2n) \\
\end{array}\right)$$
Consequently, if $T= \sum u_k^* \otimes v_k$ is a nuclear operator on $Z_2$ one gets
$$T = \left(\begin{array}{cc}
\sum_k \mathbf u_k(2n-1)^* \otimes \mathbf v_k(2n-1)& \sum_k \mathbf u_k(2n)^*\otimes \mathbf v_k(2n-1) \\
\sum_k \mathbf u_k(2n-1)^*\otimes \mathbf v_k(2n) & \sum_k \mathbf u_k(2n)^*\otimes \mathbf v_k(2n) \\
\end{array}\right)$$

\subsection{Operators acting on the scale of  $\ell_p$ spaces}\label{scale} The fact that $Z_2$ is the derived space at $\ell_2$ for the scale of $\ell_p$ spaces provides some elements of $\mathfrak L(Z_2)$. We say that an operator $\alpha:\ell_2\To\ell_2$ \emph{acts on the scale} when there are $1\leq p<2<q\leq \infty$ such that both $\alpha: \ell_p \to \ell_p$ and $\alpha: \ell_q\to \ell_q$ are bounded. In this case
$\tau_\alpha = \left(\begin{array}{cc}
\alpha & 0 \\
0 & \alpha \\
\end{array}\right)$
is a continuous operator on $Z_2$ (see \cite{symmetries21}).\medskip

Against the naive intuition, $\tau_\alpha$ can be a bounded operator on $Z_2$ and still $\alpha$ is not necessarily an operator acting on the scale, as the following simple example shows: Let $z\in\ell_f$ such that $z\not\in\ell_p$ for $p<2$. Then $\alpha_0(x)= e^*_1(x) z$ defines a bounded operator on $\ell_2$ but $\alpha_0(\ell_p)\not\subset\ell_p$ for $p<2$; thus $\alpha_0$ does not act in the scale. However $\tau_{\alpha_0}$ is a bounded operator on $Z_2$ because $(0,e^*_1), (e^*_1,0)\in Z_2^*$, $(z,0), (0,z)\in Z_2$ and $\tau_{\alpha_0}= (0,e^*_1)\otimes (z,0) + (e^*_1,0)\otimes (0,z)$.\medskip

We next present several natural examples of operators $\alpha$ acting on the scale:
 \begin{enumerate}
\item $\alpha$ a diagonal operator $D_\sigma$ with $\sigma\in \ell_\infty$, or $\alpha$ a right (or left) shift operator.
\smallskip

\item $\alpha$ a surjective isometry on $\ell_p$ for some $p\neq 2$.
\quad
It was proved in \cite[Proposition 2.f.14]{lindtzaf} that these operators have the form  $\alpha\big((x_n)_n\big)= (\varepsilon_n x_{\sigma(n)})_n$ for some permutation $\sigma:\mathbb{N} \rightarrow\mathbb{N}$ and some sequence of signs $(\varepsilon_n)_n$.
The induced operator $\tau_\alpha$ is then an isometry \cite{kaltdiff}.
\smallskip

\item $\alpha$ the \emph{Ces\`{a}ro operator} $C$ defined by $C\big((x_n)_n\big)= \Big(\frac{1}{n}\sum^n_{k=1}x_k\Big)_n$.

It is bounded on $\ell_p$ for $p>1$ with $\|C\|_p=\frac{p}{p-1}$ \cite{Rhoades_1971}. It is not bounded on $\ell_1$ since $C(e_1)$ is the harmonic series. See \cite{Brown_Halmos_Shields} for the properties of $C$ as an operator on $\ell_2$.
\smallskip

\item $\alpha$ a \emph{Hilbert matrix operator} $H_\lambda$ defined by a Hilbert matrix
$$
\Big(\frac{1}{n+m+\lambda} \Big)_{n,m=0}^\infty,\quad \lambda\in\C\setminus \{0,-1,-2,\ldots\}.
$$
The operator $H_\lambda$ is bounded on $\ell_p$ for $1<p<\infty$. See \cite{Silbermann}.
\smallskip

\item $\alpha$ a \emph{Hausdorff operator} $A$  associated to a sequence of complex scalars $\{\mu_n : n=0,1,2,\ldots\}$ (see \cite[Section 3.4]{Boos:00}). It has the form $A\big((x_k)_k\big)= \Big(\sum_{k=0}^n a_{nk}x_k\Big)_n$ with
$$a_{nk}=\begin{cases}
\binom{n}{k}\Delta^{n-k}\mu_k & 0\leq k\leq n \\
0 & k>n \end{cases},$$
where $\Delta(\mu_n)= \mu_n-\mu_{n+1}$. Many Hausdorff operators are bounded on some spaces $\ell_p$ \cite{Boos:00, Rhoades_1971}:
\begin{itemize}
\item[(i)]{the \emph{generalized Ces\`{a}ro operator} $C_{a}^\alpha$ arising from  $\mu_n=\frac{\Gamma(a+\alpha)\,\Gamma(n+a)}{\Gamma(a)\,\Gamma(n+a+\alpha)}$, where $\Gamma$ is  the Gamma function. For $\alpha>0$, $a>1/p$ the operator $C_{a}^\alpha$ is bounded on  $\ell_p$ ($p>1$) with norm equal to $\frac{\Gamma(a+\alpha)\, \Gamma(a-1/p)}{\Gamma(a+\alpha-1/p)}$.} Note that $C_{1}^1=C$, the Ces\`aro operator.
\item[(ii)] {the \emph{H\"{o}lder operator} $H_\alpha$ arising from $\mu_n=(n+1)^{-\alpha}$. For $\alpha>0$, $H_\alpha$ is bounded on $\ell_p$ with norm ${(\frac{p}{p-1})}^\alpha$;}
\item[(iii)]{the \emph{Euler operator} $(E,r)$ arising  from $\mu_n=a^n$, where $0<a<1$, $r=\frac{1-a}{a}$ and its norm on $\ell_p$ is $(1+r)^{1/p}$;}
\item[(iv)]{the \emph{Gamma operator} $\Gamma_a^\alpha$  arising from $\mu_n= \big(\frac{a}{n+a}\big)^\alpha$, where $\alpha>0$, and it is bounded on $\ell_p$ with norm $\big(\frac{a}{a-1/p}\big)^\alpha$ whenever $a>1/p$.}
\end{itemize}
\end{enumerate}

We do not know whether every Hausdorff operator which is bounded on $\ell_2$ acts on the scale.
\medskip

The following result belongs to Sneiberg \cite{Sneiberg}: \emph{Let $(X_0,X_1)$ and $(Y_0,Y_1)$ be two interpolation pairs such that $T:X_i\rightarrow Y_i$ is bounded for $i=0,1$. If $T^{-1}:X_\theta \rightarrow Y_\theta$ exists and it is bounded for some $0<\theta<1$, then there is $\varepsilon>0$ such that $T^{-1}:X_s\rightarrow Y_s$ exists and it is bounded for $|s-\theta|<\varepsilon$}. We can infer from that:

\begin{prop}\label{T_04}
For an operator $\alpha:\ell_2\to\ell_2$ acting on the scale, the following statements are equivalent:
\begin{itemize}
\item[(i)] $\tau_\alpha$ is an isomorphism.
\item[(ii)] $\alpha:\ell_2\rightarrow\ell_2$ is an isomorphism.
\item[(iii)] There exists $\varepsilon>0$ such that $\alpha:\ell_p\rightarrow\ell_p$ is an isomorphism for all $|2-p|<\varepsilon$.
\end{itemize}
Consequently $\sigma(\tau_\alpha)=\sigma(\alpha)$.
\end{prop}
\begin{proof} If $\alpha$ acts on the couple $(\ell_p,\ell_{p^*})$ then the operator $\alpha -\lambda I$ also acts on that same couple.
Moreover, if $\tau_\alpha$ is an isomorphism then $\alpha: \ell_2\to \ell_2$ is an isomorphism. And that if \emph{both} $\alpha$ and $\alpha^{-1}$ act on some scale $(\ell_p,\ell_{p^*})$ then $\tau_\alpha$ is an isomorphism on $Z_2$.\end{proof}

If $\alpha$ is an operator acting on the scale, its spectrum on $\ell_p$ maybe independent of $p$, as it is the case of diagonal operators or the left and right shift operators, but it also may vary with $p$: Leibowitz \cite{Leibowitz_72} proved that, for $1<p<\infty$ the Ces\`{a}ro operator $C$ on $\ell_p$ has no eigenvalues and its spectrum is
$$
\sigma(C)=\big\{\lambda\in\mathbb{C}\colon \big|\lambda-p^*/2\big|\leq p^*/2\big\}.
$$
Moreover, $\lambda I - C$ is a Fredholm operator with index $-1$ for $|\lambda-p^*/2|< p^*/2$, and has dense proper range for $|\lambda-p^*/2|= p^*/2$ \cite{Cesaro-85}.\smallskip

The Hilbert matrix operator $H_1$ on $\ell_2$ has no eigenvalues and $\sigma(H_1)$ is the interval $[0,\pi]$, see \cite{Magnus}, while the spectrum on $\ell_p$ varies with $p$, and it has eigenvalues for $p>2$ and residual points for $p<2$ \cite{Silbermann}. Since the matrix representing $H_1$ is symmetric, the conjugate operator of $H_1:\ell_p\to \ell_p$ is $H_1:\ell_{p^*}\to \ell_{p^*}$. Thus the spectra of $H_1$ on $\ell_p$ and $\ell_{p^*}$ coincide.


\subsection{Operators on the Calder\'on space}\label{opcal}
We obtain operators on $Z_2$ by picking operators on the Calder\'on space $T: \mathcal C\to \mathcal C$ such that $T[\ker \delta \cap \ker \delta' ] \subset \ker \delta \cap \ker \delta'$. The simplest way to do that is to pick an operator $\tau$ on the scale and then set $T(f)(z)=\tau(f(z))$. If $\varphi: \mathbb S\to \mathbb D$ is a conformal map, then the operator $S(f)(z) = \tau \left( \varphi(z) f(z)\right)$ induces
$\left(\begin{array}{cc}
0 & \tau \\
0 & 0 \\
\end{array}\right).$
Therefore, given two operators $\alpha, \phi$ on the scale, $T(f)(z) = \alpha(f(z)) + \beta(\varphi(z)f(z))$ induces the upper triangular operator
$\left(\begin{array}{cc}
\alpha & \beta \\
0 & \alpha \\
\end{array}\right)$ on $Z_2$.

\subsection{Diagonal operators} The continuity of diagonal operators $D_\sigma$ on $Z_2$ is  related with the unconditional structure of the space. Recall that the sequence $(u_n)_{n\in\mathbb{N}}$ given by $u_{2n-1}=(e_n,0)$ and $u_{2n}=(0,e_n)$, where $(e_n)_n$ is the canonical basis of $\ell_2$, is a basis on $Z_2$ which is not unconditional. Therefore not all $\sigma \in\ell_\infty$ define a diagonal operator on $Z_2$. Let us denote $a=(\sigma_{2n -1})$ and $b= (b_n=\sigma_{2n})$.
If $D_\sigma = \left(\begin{array}{cc}
D_a & 0 \\
0 & D_b \\
\end{array} \right)$
is an operator on $Z_2$ then $D_a - D_b= D_{a-b}$ is compact by Lemma \ref{ncfutm}; thus $a-b\in c_0$. It is startling that an additional condition is required:
\begin{prop} A diagonal operator $D_\sigma: Z_2\rightarrow Z_2$ defined by a monotone decreasing sequence $\sigma$ is bounded if and only if  $\big(|\sigma_{2k-1}-\sigma_{2k}|\big)_{k\in\mathbb{N}}= \mathrm{O}\Big(\frac{1}{\log n}\Big)$.\end{prop}
\begin{proof} Cabello and Garc\'ia showed in \cite{cabegar} that a diagonal operator $D_a :\ell_2\rightarrow\ell_2$ can be lifted to $Z_2$ if and only if the decreasing rearrangement sequence $(a_n^*)_n$ is $\mathrm{O}\Big(\frac{1}{\log n}\Big)$. The self duality of the Kalton-Peck space yields the result. \end{proof}

\subsection{Block operators}\label{blocks} Let $U=(u_n)$ be a bounded sequence of disjointly supported blocks in $\ell_2$. We define a bounded operator $u:\ell_2\to \ell_2$ by $u x=\sum x_n v_n$.
Kalton \cite{kaltsym} defined the operator $T_U$ on $Z_2$ by $T_U(e_n,0)= u_n$ and $T_U(0,e_n)= (\KP u_n, u_n)$, and proved it is and into isometry. Let us call $T_U$ a \emph{block operator.} As we said in \ref{scale}, if $\alpha$ is an operator acting on the scale then  $\tau_\alpha$ is an upper triangular operator on $Z_2$.
In general, a perturbation of $\tau_\alpha$ is required to make it an upper triangular operator. In particular, the operator $u$ defined by a sequence of disjointly supported normalized blocks in $\ell_1$ is not an operator on the scale.  In \cite[Theorem 4.6]{gspaces} it is explained how to obtain the required perturbation and how this perturbation yields the Kalton block operator
$T_U = \left(\begin{array}{cc}
u & \KP u \\
0 & u \\
\end{array}\right)$
mentioned above.

\section{Operator ideals on $Z_2$} \label{sect:ideals}

The classes $\mathfrak S$ and $\mathfrak C$ are not  dual to each other. In general $T^*\in \mathfrak S \Longrightarrow T\in \mathfrak C$ and $T^*\in \mathfrak C \Longrightarrow T\in \mathfrak S$. But since $Z_2$ is reflexive, it turns out that an operator $T:Z_2\rightarrow Z_2$ is strictly singular (resp. cosingular) if and only if $T^*:Z_2^* \rightarrow Z_2^*$ is strictly cosingular (resp. singular). One moreover has:

\begin{theorem} One has the identities
\label{singcosing}
$$\mathfrak S(Z_2) = \mathfrak C(Z_2) =\mathfrak{In}(Z_2).$$
Moreover, that set contains every proper ideal of $\mathfrak L(Z_2)$.
\end{theorem}
\begin{proof}
By the first part of Theorem \ref{Kalton}, if $S\in \mathfrak L(Z_2)\setminus \mathfrak S(Z_2)$ then there exists $A,B\in \mathfrak L(Z_2)$ so that $ASB=I_{Z_2}$; hence $S$ does not belong to any proper operator ideal, and $\mathfrak{In}(Z_2) \subset  \mathfrak S(Z_2)$. Since $\mathfrak S(X)$ and $\mathfrak C(X)$ are contained in $\mathfrak{In}(X)$ for each $X$ and $Z_2\simeq Z_2^*$, the equalities follow.
%
%
\end{proof}

Observe that $\mathfrak S(Z_2,\ell_\infty)) \neq \mathfrak L(Z_2,\ell_\infty)= \mathfrak{In} (Z_2,\ell_\infty)$: let $T\in\mathfrak L(Z_2,\ell_\infty)$. For every $A\in \mathfrak L(\ell_\infty,Z_2)= \mathfrak S (\ell_\infty,Z_2)$, $I-AT$ is Fredholm.
Similarly,  $\mathfrak C(\ell_1,Z_2)) \neq \mathfrak L(\ell_1,Z_2)= \mathfrak{In} (\ell_1,Z_2)$.\medskip

The next result is a dual version of the second part of Theorem \ref{Kalton}.

\begin{prop}\label{Kalton-dual}
Let $T\in \mathfrak L(X,Z_2)$.
If $T\notin \mathfrak C$ then there exists a complemented subspace $N$ of $Z_2$ with $Z_2/N$ isomorphic to $Z_2$ such that $Q_NT$ is surjective, where $Q_N:Z_2\to Z_2/N$ is the quotient map.
\end{prop}
\begin{proof}
If $T\in \mathfrak L(X,Z_2)$ is not in $\mathfrak C$ then $T^*\in \mathfrak L(Z_2^*, X^*)$ is not in $\mathfrak S$. Since $Z_2^*\simeq Z_2$, by Theorem \ref{Kalton} there exists a complemented subspace $M$ of $Z_2^*$ isomorphic to $Z_2^*$ such that $T^*|_M$ is an isomorphism. Then $N=\, ^{\perp} M$ is a subspace of $Z_2$ satisfying the required conditions.
\end{proof}

Let $\mathfrak K$ be the class of compact operators and let $L_p\equiv L_p(0,1)$ for $1\leq p\leq \infty$. Then $\mathfrak S(L_p)\neq \mathfrak K(L_p)$ for $p\neq 2$ \cite{Milman:70}, but $T\in \mathfrak S(L_p)$ implies $T^2 \in \mathfrak K(L_p)$ \cite{GMF:60}.

\begin{theorem}
We have $\mathfrak S(Z_2)\neq \mathfrak K(Z_2)$, but $S,T\in \mathfrak S(Z_2)$ implies $ST \in \mathfrak K$.
\end{theorem}
\begin{proof}
As we mentioned before, $i\, p\in \mathfrak{L} (Z_2)$ is strictly singular but not compact.

For the remaining part, recall that an operator $S$ acting on a reflexive space $X$ is compact if and only if for every normalized weakly null sequence $(x_n)$ in $X$, $(Sx_n)$ has a norm null subsequence; and it was proved in \cite[Theorem 5.4]{kaltpeck} that every normalized weakly null sequence $(x_n)$ in $Z_2$ has a subsequence equivalent either to the (usual) basis of $\ell_2$ or to the (usual) basis on $\ell_f$.

Let $S,T\in \mathfrak S(Z_2)$ and let $(x_n)$ be a normalized weakly null sequence $Z_2$. If $(x_{n_k})$ is a subsequence equivalent to the basis of $\ell_f$ then $(Tx_{n_k})$ has no  subsequence equivalent to the basis of $\ell_f$ because $T$ is strictly singular; hence $(Tx_{n_k})$ has a   subsequence equivalent to the basis of $\ell_2$ or it is norm null. Also, if $(x_{n_k})$ is a subsequence equivalent to the basis of $\ell_2$ then $(Tx_{n_k})$ has no  subsequence equivalent to the basis of $\ell_2$, and has no  subsequence equivalent to the basis of $\ell_f$ because $\ell_f \subsetneq \ell_2$: a bounded operator cannot take the unit basis of $\ell_2$ to the unit basis of $\ell_f$; hence it is norm null. In each case,  $(STx_n)$ has a norm null subsequence.
\end{proof}

The \emph{perturbation class} of a class of operators $\mathcal{A} \subset \mathfrak L$ is defined by its components as follows when $\mathcal A(X,Y)\neq\emptyset$: $$P\mathcal A(X,Y)= \{L\in \mathfrak L(X,Y)\colon T+L\in \mathcal{A}(X,Y)\text{ for all }T\in \mathcal{A}(X,Y)\}.$$

Kato and Vladimirskii (see \cite[Section 26.6]{pietsch}) proved that $\mathfrak S \subset P \Phi_+$ and $\mathfrak C\subset P \Phi_-$, and it is known that $\mathfrak{In} = P\Phi$. The \emph{perturbation classes problem} asks whether $\mathfrak S=P\Phi_+$ and $\mathfrak C=P\Phi_-$. This  problem has a positive answer under certain conditions but not in general (see \cite{Gon_pert}), and also for $Z_2$ although this space does not verify  those conditions.

\begin{prop} We have $P\Phi(Z_2) =P\Phi_+(Z_2) =P\Phi_-(Z_2) = \mathfrak S (Z_2)$.
\end{prop}
\begin{proof} In general, $\mathfrak S(X)\subset P\Phi_+(X)\subset P\Phi(X) =\mathfrak{In}(X)$ and $\mathfrak C(X)\subset P\Phi_-(X) \subset \mathfrak{In}(X)$,
but Theorem \ref{singcosing} implies $\mathfrak S(Z_2)= \mathfrak C(Z_2)= \mathfrak{In}(Z_2)$. So all these  inclusions are equalities for $X=Z_2$.
\end{proof}


\section{Further directions of research}
\label{sect:problems}

The overall tone of these suggestions is to determine which properties of operators on $\ell_2$ are valid for operators on $Z_2$ and which are not.

\subsection{The convolution on $\mathfrak{L}(Z_2)$}

Here we consider the relation between $T$ and $T^+$ as operators in $\mathfrak{L}(\ell_2)$.

\begin{quest}\label{T^+T}
Suppose that $T\in \mathfrak{L}(Z_2)$ is an isomorphism into. Is $T^+T$ bijective? What does it mean $T^+T=I$ or $TT^+T=T$ for $T\in \mathfrak{L}(Z_2)$?
\end{quest}


\begin{quest}\label{B*-algebra}
Is $\mathfrak{L}(Z_2)/\mathfrak S(Z_2)$ isomorphic to a C$^*$-algebra?
\end{quest}

Clearly $\mathfrak{L}(Z_2)$ is not isomorphic to a C$^*$-algebra since there exists $T\neq 0$ such that $T^+T=0$. However $T^+T\in \mathfrak S(Z_2)$ implies $T\in \mathfrak S(Z_2)$, and $T\in \mathfrak S(Z_2)$ if and only if $T^+\in \mathfrak S(Z_2)$.

\subsection{Polynomially bounded operators}
A \emph{contraction} is an operator $T$ with $\|T\|\leq 1$, and for a polynomial $p$ we denote $\|p\|_\infty=\sup_{|z|<1} |p(z)|$.
\begin{quest}\label{polyn-bounded}
Is every contraction in $\mathfrak{L}(Z_2)$ polynomially bounded? Equivalently,
\begin{equation}\label{vN-ineq}
\exists\, C>0 \textrm{ so that } \|T\|\leq 1\Rightarrow \|p(T)\|\leq C\|p\|_\infty \textrm{ for every polynomial }p?
\end{equation}
\end{quest}
Note that (\ref{vN-ineq}) with $C=1$ isometrically characterizes Hilbert spaces. Moreover, if $X$ is isomorphic to a Hilbert space clearly (\ref{vN-ineq}) holds for some $C\geq 1$; however the converse implication fails (see \cite{Zarrabi}).

\subsection{The group of invertible operators}
We denote by $\mathfrak{GL}(X)$ the group of invertible operators on a Banach space $X$. It is known that $\mathfrak{GL}(\ell_2)$ is connected in the complex case, while $\mathfrak{GL} (\ell_p \times\ell_q)$ is not connected for $1\leq p<q<\infty$ \cite{em,neu}.

\begin{quest}\label{GL(Z2)}
In the case $\mathbb K=\C$, is $\mathfrak{GL}(Z_2)$ connected?

Is the subgroup $\{T\in \mathfrak{GL}(Z_2) : T \textrm{ is upper triangular}\}$ connected?
\end{quest}

The latter question could be tackled by obtaining a characterization of the invertible operators $T\in \mathfrak L(Z_2)$ in terms of the components $\alpha, \beta, \delta, \gamma$ of the matrix representation of $T$.

\subsection{Representations of $Z_2$} A basic question whose meaning is not even clear is whether there are other ``natural" presentations of $Z_2$ beyond the $\ell_2$ and the $\ell_f$ presentations considered in this paper. In homological terms, since $Z_2 \simeq Z_2 \oplus Z_2$ one could obtain other nontrivial representations such as $\xymatrix{0 \ar[r]& \ell_2 \ar[r]& Z_2 \ar[r] & \ell_2 \oplus Z_2 \ar[r]&0},$
$\xymatrix{0 \ar[r]& \ell_f \ar[r]& Z_2 \ar[r] & \ell_f^* \oplus Z_2 \ar[r]&0},$
etc. Or even the trivial one
$$\xymatrix{0 \ar[r]& Z_2 \ar[r]& Z_2 \ar[r] & Z_2 \ar[r]&0}.$$
None of these representations are even ``isomorphic" to either $(\ref{ip-seq})$ or $(\ref{jq-seq})$.




\end{document}